\documentclass[12pt,twoside]{article}
\usepackage{a4,latexsym}

\font\Bbbb=msbm7
\newcommand{\bZ}{\mathbf{Z}}
\newcommand{\bA}{\mathbf{A}}

\newtheorem{theo}{Theorem}[section]

\newtheorem{lema}{Lemma}[section]
\newtheorem{rema}{Remark}[section]

\def \al{\alpha}  \def \be{\beta}    \def \ga{\gamma}
   \def \eps{\varepsilon}

 \def \Ga{\Gamma}   \def\vt{\vartheta}

\renewcommand {\div}{{\rm div}\,}

\newcommand {\ee}{{\rm e}\,}
\newcommand {\sign}{{\rm sign}\,}
 1
 2
\font\Bbbb=msbm10 at 12pt
\font\Bbbbb=msbm7 scaled \magstep1

 3
 5

\newcommand {\rtt}{\mbox{\Bbbbb R}}

\newcommand {\rdd}{\mbox{\Bbbb R}}

\newcommand {\bnul}{{\mathbf {0}}}
\newcommand {\bv}{{\mathbf{v}}}

\newcommand {\bu}{{\mathbf{u}}}
\newcommand {\bI}{{\mathbf{I}}}
\newcommand {\bM}{{\mathbf{M}}}

\newcommand {\bff}{{\mathbf{f}}}
\newcommand {\bx}{{\mathbf{x}}}
\newcommand {\bm}{{\mathbf{m}}}

\newcommand {\bb}{{\mathbf{b}}}

\newcommand {\bp}{{\mathbf{p}}}

\newcommand {\bQ}{{\mathbf{Q}}}

\newcommand {\bF}{{\mathbf{F}}}

\newcommand {\pder}[2]{\frac {\partial #1}{\partial #2}}
\newcommand {\demo}{{\bf Proof:} }
\def\kondemo{\hfill \vrule width7pt height7pt depth0pt}

\font\index = cmr7
\newcommand{\wor}[2]{W^{#1,#2}(\rdd^3)}

\renewcommand{\lor}[1]{L^{#1}(\rdd^3)}

\begin{document}

\title{A counterexample to the smoothness of the solution to an equation
arising in fluid mechanics 
\footnotetext{
\hfil\break\indent
}}

\author {Stephen Montgomery-Smith\thanks{Partially supported by 
the National Science Foundation DMS 9870026, and a grant from the
Research Board of the University of Missouri.}
 \ \ Milan Pokorn\'y\thanks{Supported by the Grant Agency of the 
Czech Republic (grant
No. 201/00/0768) and by the Council of the Czech Government
(project No.~113200007)}
\date{}
}
\maketitle

\begin{abstract}
We analyze the equation coming from the Eulerian-Lagrangian description
of fluids.  We discuss a couple of ways to extend this notion to 
viscous fluids.  The main focus of this paper is to discuss 
the first way, due to Constantin.  We show that this description
can only work for short times, after which the ``back to coordinates map''
may have no smooth inverse.  Then we briefly discuss a second way that
uses Brownian motion.  We use this to provide a plausibility argument
for the global regularity for the Navier-Stokes equations.
\end{abstract}

\section{Introduction}

Recently there has been interest in some new variables describing 
the solutions to the Navier-Stokes and Euler equations.  These variables
go under various names, for example, {\em the magnetization variables},
{\em 
impulse variables}, {\em velicity} or {\em Kuzmin-Oseledets variables}.

Let us start by considering the incompressible Euler equations in the entire 
three-dimensional space, that is,
\begin{equation} \label{euleru} %1.1
\begin{array}{c}
\left. \begin{array}{c}
\displaystyle \pder{\bu}{t} + \bu\cdot \nabla \bu + \nabla p = \bff \\[6pt]
\div \bu = 0 
\end{array} \right\} {\rm in \ } \rdd^3 \times (0,T) \\[10pt]
\bu(\bx,0) = \bu_0(\bx) \quad {\rm in \ } \rdd^3\, ,
\end{array}
\end{equation}
where $\bu$ and $p$ are given functions, $\bu: \rdd^3 \times(0,T) 
\mapsto \rdd^3$ and $p: \rdd^3 \times(0,T) 
\mapsto \rdd$, $0<T\leq \infty$ (see Section 1 for further 
explanation).

The question of global existence of even only weak solutions to 
system~(\ref{euleru}) is an open question and only the existence of either 
measure-valued solutions (see \cite{DPMa}) or dissipative solutions (see
\cite{Lio}) is known. Nevertheless, a common approach to try to prove 
the global existence of smooth solutions is to use local existence results,
and thus reduce the problem to proving a priori estimates.  So we will
assume that we have a smooth solution to the equations.

In that case, we can rewrite the Euler equations as the following system
of equations (see for example \cite{Ch}):
\begin{equation} \label{eulerm} 
\begin{array}{c}
\left. \begin{array}{c}
\displaystyle \pder{\bm}{t} + \bu\cdot \nabla \bm + \bm \cdot (\nabla \bu)^T = 
  \bff \\[6pt]
\bu = \bm - \nabla \eta \\[6pt]
\div \bu = 0 
\end{array} \right\} {\rm in \ } \rdd^3 \times (0,T) \\[10pt]
\bm(\bx,0) = \bu_0(\bx) \quad {\rm in \ } \rdd^3\, .
\end{array}
\end{equation}
Here $\bm:\rdd^3 \times(0,T) 
\mapsto \rdd^3$ is called the magnetization variable.  This new formulation
has several advantages to the usual one, in particular the solution
can be written rather nicely in the following way.  
Suppose that
the initial value for $\bm$ may be written as
\begin{equation}
\label{init-m-alpha-beta}
\bm(\bx,0) = \sum_{i=1}^R \beta_i(\bx,0) \nabla \alpha_i(\bx,0) ,
\end{equation}
and suppose that $\alpha$ and $\beta$ satisfy the transport equations,
that is
%\begin{equation} \label{transpalphabeta}
$$
\begin{array}{c}
\displaystyle \pder{\alpha_i}{t} + \bu \cdot \nabla \alpha_i = 0
\\[6pt]
\displaystyle \pder{\beta_i}{t} + \bu \cdot \nabla \beta_i 
= \sum_{j=1}^R Q_{j,i} f_j
\end{array}
%\end{equation}
$$
where $\bQ = (Q_{i,j})$ is the matrix inverse of the matrix whose entries are
$\pder{\alpha_i}{x_j}$.  (There is some difficulty to suppose that this
inverse exists unless $R=3$ --- see below.  But generally this will not 
be a problem if $\bff = 0$.)
Then 
%\begin{equation}
%\label{m-alpha-beta}
$$
\bm(\bx,t) = \sum_{i=1}^R \beta_i(\bx,t) \nabla \alpha_i(\bx,t) 
$$
%\end{equation}
is the solution to system~(\ref{eulerm}).  That is to say, at least in
the case that $\bff = 0$, the 
magnetization variable may be thought of as a ``1-form'' acting 
naturally under a change of basis induced by the flow of the fluid.

The advantage of the magnetization variable is that it is local in that
its support never gets larger, it is simply pushed around by the flow.
It is only at the end, after one has calculated the final value of
$\bm$, that one needs to take the Leray projection to compute
the velocity field $\bu$.

Indeed one very explicit way to write $\bm$ according to 
equation~(\ref{init-m-alpha-beta}) is to set $\alpha_i(\bx,0)$
equal to the $i$th unit vector, and $\beta_i(\bx,0) = u_i(\bx,0)$, 
for $1 \le i \le R = 3$.  In that case
let us denote $A^E_i(\bx,t) = \alpha_i(\bx,t)$ and 
$v_i(\bx,t) = \beta_i(\bx,t)$.  In that case we see that
$\bA^E(\bx,t)$ is actually the
back to coordinates map, that is, it denotes the initial position of the
particle of fluid that is at $\bx$ at time $t$ (see for example \cite{Ch}).  
Furthermore
in the case that $\bff = 0$, we see that $\bv(\bx,t) = \bu_0(\bA^E(\bx))$.
Furthermore, it is well known if $\bu$ is smooth, that
$\bA^E(\cdot,t)$ is smoothly invertible, and that the determinant of the 
Jacobian of $\bA^E$ is identically equal to $1$ (because $\div \bu = 0$).  
Hence the matrix $\bQ$ exists.  For definiteness, we
write the explicit equation for
$\bm$:
\begin{equation} \label{m-A-v}%1.2
m_i(\bx,t) = \pder{\bA^E(\bx, t)}{x_i} \cdot \bv^E(\bx, t) .
\end{equation}

The desire, then, is to try to extend this notion to the Navier-Stokes
equations
\begin{equation} \label{navier-stokes-u}%{1.3} 
\begin{array}{c}
\left. \begin{array}{c}
\displaystyle \pder{\bu}{t} + \bu\cdot \nabla \bu - \nu \Delta \bu + \nabla p = \bff \\[6pt]
\div \bu = 0 
\end{array} \right\} {\rm in \ } \rdd^3 \times (0,T) \\[10pt]
\bm(\bx,0) = \bu_0(\bx) \quad {\rm in \ } \rdd^3\, .
\end{array}
\end{equation}
(Only local-in-time existence of smooth solutions to the Navier-Stokes 
equations is known --- see for example \cite{KiLa}; 
globally in time, only  existence 
of weak solutions is known, see \cite{Le}.)

Again, these can be rewritten into the magnetic variables formulation
as follows:
\begin{equation} \label{navier-stokes-m}%{1.3a} 
\begin{array}{c}
\left. \begin{array}{c}
\displaystyle \pder{\bm}{t} -\nu \Delta \bm + \bu\cdot \nabla \bm + \bm \cdot (\nabla \bu)^T = 
  \bff \\[6pt]
\bu = \bm - \nabla \eta \\[6pt]
\div \bu = 0 
\end{array} \right\} {\rm in \ } \rdd^3 \times (0,T) \\[10pt]
\bu(\bx,0) = \bu_0(\bx) \quad {\rm in \ } \rdd^3\, .
\end{array}
\end{equation}
The problem is to find the analogue of equation~(\ref{m-A-v}).  The difficulty
comes from the term $\nu \Delta \bm$.  There are two ways known to the
authors --- one is to use probabilistic techniques.  Since this technique
seems to be not as widely known as it should be, we will include a 
short (non-rigorous) description of this method at the end of the paper.
We will also include a short plausibility argument for the global regularity
for the Navier-Stokes equations.

Another approach was developed by Peter Constantin (see \cite{Co}).  He
used new quantities $\bA^N$ and $\bv^N$ obeying the following equations.
Let us represent $\bu$ in a form similar to 
(\ref{eulerm})
%\begin{equation} \label{1.4}
$$
u_i(\bx,t) = \pder{\bA^N(\bx,t)}{x_i} \cdot \bv^N(\bx,t) - \pder{n(\bx,t)}{x_i}\, ,
%\end{equation}
$$
where 
\begin{equation} \label {0.3}
\begin{array}{c}
\Ga(\bA^N) = \bnul \qquad {\rm in \ } \rdd^3 \times (0,T) \\[6pt]
%\end{equation}
%\begin{equation} \label{0.4}
\bA^N(\bx, 0) = \bx \qquad {\rm in \ } \rdd^3 
\end{array}
\end{equation}
$$
\Ga = \pder{}{t} + \bu\cdot \nabla - \nu \Delta\, ,
$$
and $\bv^N$ obeys a rather complicated equation
$$
\Gamma (v^N_i) = 2 \nu C_{m,k;i} \pder{v_m}{x_k} + Q_{j,i} f_j\, ,
$$
where $\bQ$ is the inverse matrix to $\nabla \bA^N$ and $\Gamma_{i,j}^m =
- Q_{k,j} C_{m,k;i}$ denotes the Christoffel coefficients.  
In order for the equation
for $\bv$ to make sense, it is necessary for the map $\bA^N$ to have a
smooth inverse.  An approach to proving such a result is to consider
the system of PDE's
\begin{equation} \label{Q-equ} % {0.1}
\begin{array}{c}
\Ga (\bQ) = (\nabla \bu) \bQ + 2\nu \bQ \partial_k (\nabla \bA^N) \partial_k \bQ
\qquad {\rm in \ } \rdd^3 \times (0,T) \\[6pt]
%\end{equation}
%\begin{equation}  \label{Q-init} % {0.2}
\bQ(\bx,0) = \bI \qquad {\rm in \ } \rdd^3\, ,
\end{array}
\end{equation}
If the above equations have smooth solutions, then it can easily be shown
that $\bQ$ is the inverse to $\nabla \bA^N$.  However, the problem is that
while it is easy and standard to show that equation~(\ref{Q-equ}) has local
smooth solutions, it is not clear that it has global solutions in any
sense at all.

The purpose of this note is to show that indeed global smooth solutions
do not exist.  As Peter Constantin pointed out to us, this does not invalidate
his method, but it does mean that to make his method work for a large
time period that one has to break that
interval into shorter pieces, and apply the method to each small interval.

The main result is summarized in the following theorem.
\begin{theo} \label {t 0.1}
There exists $\bu \in C^\infty_0(\rdd^3 \times [0,\infty))$, divergence free such 
that if $\bA^N$ is a smooth solution to (\ref{0.3})
%--(\ref{0.4}) 
then there exists $t>0$ such that
\begin{itemize}
\item [(a)] $\bA^N(\bnul, t)$ does not have a smooth inverse
\item [(b)] $\limsup_{\tau \to t^-} \|\bQ\|_\infty(\tau) = \infty$,
\end{itemize}
where $\bQ$ is a solution to (\ref{Q-equ})
%--(\ref{Q-init}) 
corresponding to $\bu$ 
and $\bA^N$.
\end{theo}

\section{Outline of the Proof of Theorem \ref{t 0.1}}   
\setcounter{equation}{0}

In this section we will give the plan for the proof of Theorem~\ref{t 0.1}.
The idea of the proof is really quite simple.  We will in fact construct
a family of divergence free, smooth solutions $\bu_s$ 
to (\ref{Q-equ})
%--(\ref{Q-init})
parameterized by a number $s \in [0,2\pi]$.  We will use simple ideas
from algebraic topology to show that there exists $s_0 \in [0,2\pi]$ such that
$u_{s_0}$ provides an example to prove Theorem~\ref{t 0.1}.

In fact all of the solutions we construct will be axisymmetric, indeed,
when written in cylindrical coordinates, they have the form:
$\bu_s = (0, u_\theta(r,z,t), 0)$.
We will prove that there exists $s_0 \in [0,2\pi]$ and $t_0>0$ 
such that the Jacobian $\nabla \bA^N_{s_0}(0,t_0)$ is non-invertible.
To this end we have the following representation result.

\begin{lema} \label{c 5.1}
For any $t>0$, $\nabla \bA^N_s(\bnul, t)$ can be 
uniquely written as
%\begin{equation} \label{5.1}
$$
\left(\begin{array}{lll}
a_s \cos b_s, & -a_s \sin b_s & 0\\
a_s \sin b_s, &  a_s \cos b_s & 0\\
0 & 0 & 1
\end{array} \right)
%\end{equation}
$$
for some $a_s(t)\in \rdd^+$ and $b_s(t) \in [0,2\pi)$.
\end{lema}

The proof of Theorem~\ref{t 0.1} will proceed as follows.
For each $s \in [0,2\pi]$, we will construct $\bu_s$.  The Theorem
will be proved if we can show the existence of $s_0 \in [0,2\pi]$
and $t_0 > 0$ such that $a_{s_0}(t_0) = 0$.  We will assume the opposite,
and give a proof by contradiction.

We will need some simple facts from algebraic
topology.  We refer the reader to \cite{Mau} for more details.
Let us consider the collection of continuous functions 
$[0,\infty] \to \rdd^2-\{(0,0)\}$ which map $0$ and $\infty$ to $(1,0)$.
We will say two such functions $f$ and $g$ 
are \emph{homotopic with base point $(1,0)$} (or simply \emph{homotopic}) 
if there exists a jointly continuous 
function $F:[0,\infty] \times [0,2\pi] \to \rdd^2-\{(0,0)\}$ such that
$F(\cdot,0) = f$, $F(\cdot,2\pi) = g$ and
$F(0,\cdot) = F(\infty,\cdot) = (1,0)$.  We will call the function $F$
a homotopy.  Clearly being homotopic is an
equivalence relation.  

It is well known that a constant map
$f(t) = (1,0)$, and a map with ``winding number 1'', for example,
$g(t) = (\cos(2\pi t/(1+t)),\sin(2\pi t/(1+t)))$ are not homotopic.
(Since $\rdd^2-\{(0,0)\}$ is homotopy equivalent to the unit circle,
this is basically saying that the fundamental group of the unit circle
is non-trivial.)

In order to provide our contradiction we will prove the following result.

\begin{lema}
\label{homotopy}  If $\bu_s$ is constructed as described in the
next section, with the various parameters chosen appropriately, 
then the function
$$ F(t,s) = (a_s(t) \cos b_s(t), a_s(t) \sin b_s(t)) $$
provides a homotopy between the function $f$
and a function homotopic to $g$.
\end{lema}

\section{Properties of the operator $\Ga$}
\setcounter{equation}{0}

We will not prove the smoothness of solution to (\ref{0.3});
%--(\ref{0.4}) 
it can be done in a very standard way, using the estimates to parabolic
equations given for example in \cite{LaSoUr}. 
Let us only summarize the main result
here.  This will show that the function $F$ described in 
Lemma~\ref{homotopy} is continuous on any compact subset of
$[0,\infty)\times[0,2\pi]$.

\begin{lema} \label {l 2.1}
Let $\bu \in C^\infty_0([0,T)\times \rdd^3)$ for some $T>0$. Then,
in the class of functions $V_k= \{\bv \in L^2((0,T); L^{2}_{loc}(\rdd^3)); \bv-\bx
\in L^2((0,T); W^{k,2}(\rdd^3)); \pder{\bv}{t} \in L^2((0,T);
W^{k-2,2} (\rdd^3))\}$, $k\geq 2$, there 
exists exactly one solution to (\ref{0.3}).
%--(\ref{0.4}). 
Moreover, this solution
is smooth, that is, in $C^\infty((0,T]\times \rdd^3) \cap 
C([0,T]\times \rdd^3)$, 
and $\bA^N-\bx \in L^2((0,T);
\wor{k}{2})$ for any $k\geq 0$.  Furthermore the solution depends smoothly upon
the choice of $\bu$.
\end{lema}

\begin{rema} \label{r 2.1}
{\rm Note that if $\bu$ belongs to $L^\infty((0,T); \lor{2})$ 
$\cap L^2((0,T);
\wor{1}{2})$ $\cap L^1((0,T);
\lor{\infty})$ (the usual information about a weak solution to the 
Navier-Stokes equations), 
then only $\bA^N-\bx \in L^\infty((0,T);
\wor{1}{2}) \cap L^2((0,T);
\wor{2}{2}) \cap L^\infty((0,T);
\lor{\infty})$ and $\pder{\bA^N}{t} \in L^2((0,T);
\lor{2})$. The proof is essentially the same as the proof of Lemma 
\ref{l 2.1} using \cite{LaSoUr} and is left as an exercise. }
\end{rema}

\begin{lema} \label{l 2.2}
There exists an interval $(0,t)$ such that for $\bu$ and $\bA^N$ smooth as in 
Lemma \ref{l 2.1}, $\bQ$ is a smooth solution to (\ref{Q-equ}).
%--(\ref{Q-init}).
\end{lema}

\demo The existence of the solution can be shown using the Galerkin
method combined with standard a priori estimates. We leave the details
of the proof to the reader as an exercise. 
\kondemo

\bigskip

Now, on the time interval from Lemma \ref{l 2.2} we see that
%\begin{equation} \label{2.2}
$$
\bZ = (\nabla \bA^N) \bQ - \bI
%\end{equation}
$$
obeys the equation (see \cite{Co})
\begin{equation} \label{2.3}
\Ga \bZ = 2\nu \bZ \partial_k(\nabla \bA^N)\partial_k \bQ
\end{equation}
in $\rdd^3 \times (0,T)$ with the initial condition $\bZ(\bx,0) = \bnul$. Since,
for $\nabla^2 \bA^N$ and $\nabla \bQ$ bounded, there exists the unique solution to   
(\ref{2.3}), we have $\bZ \equiv \bnul$ and thus $\bQ = (\nabla \bA^N)^{-1}$ pointwise.

\bigskip

Also, we are now in a position to prove Lemma~\ref{c 5.1}.  Since
(\ref{0.3})
%--(\ref{0.4}) 
are uniquely solvable, it follows that 
the solution is axisymmetric and hence we can apply the following result.

\begin{lema} \label{l 5.1}
Let $\bF: \rdd^2 \mapsto \rdd^2$ be a vector field which is of the 
class $C^1$ on some neighborhood of the origin and, written in polar 
coordinates, $F_r$ and $F_\vt$ are independent of $\vt$. Then
$$
\pder{F_x(\bnul)}{x} = \pder{F_y(\bnul)}{y}\, , \qquad 
\pder{F_x(\bnul)}{y} = -\pder{F_y(\bnul)}{x}\, .
$$
\end{lema}

\demo
Denote $F_r = f(r)$ and $F_\vt = g(r)$. Then we get
$$
\pder{F_x}{x} = f' \cos^2 \vt + \frac{f}{r} \sin^2 \vt - g' \sin \vt 
\cos \vt + \frac{g}{r} \sin \vt \cos \vt\, .
$$
Since $\lim_{r\to 0} \pder{F_x}{x}$ exists, necessarily 
$$
\lim_{r\to 0} \Big(f'(r) - \frac{f(r)}{r}\Big) = 0 \qquad \mbox{ and }
\qquad \lim_{r\to 0} \Big(g'(r) - \frac{g(r)}{r}\Big) = 0\, .
$$
Thus $\pder{F_x(\bnul)}{x} = f'(0)$. Next
$$
\pder{F_y}{y} = f' \sin^2 \vt + \frac{f}{r} \cos^2 \vt + g' \sin \vt 
\cos \vt - \frac{g}{r} \sin \vt \cos \vt
$$
and also $\pder{F_y(\bnul)}{y} = f'(0)$. Analogously we get that
$\pder{F_x(\bnul)}{y} = -g'(0)$ and $\pder{F_y(\bnul)}{x} = g'(0)$.
The lemma is proved. \kondemo

\section{Construction of the Fluid Flow}

\setcounter{equation}{0}
We will consider the following vector field in cylindrical coordinates:
\begin{equation}  \label{3.1}
\bu_s = \bu = (0, u_\theta(r,z,t), 0)
\end{equation}
with $u_\theta(r,z,t) = \al(r) \be(|z|) \ga_s(t) r$, where
$$
\begin{array}{l}
\begin{array}{rlllrlllll}
\al(r) &= 0 &\, {\rm for \ } & r\geq R^o, & \al(r) &= 1 &  \, {\rm for \ }
& r \leq R^i \, , & \al \in C^\infty_0([0,\infty)) \, , & 0\leq \al(r) \leq 1 \\[6pt]
\be(|z|) &= 0 &\, {\rm for \ } & |z|\geq Z^o, & \be(|z|) &= 1 &  \, {\rm for \ }
& |z| \leq Z^i \, , & \be \in C^\infty_0(\rdd) \, , & 0\leq \be(r) \leq 1 \\[6pt]
\ga_s(t) &= 0 &\, {\rm for \ } & t\geq t_0, & \ga_s(0) &= 0, & &  \ga_s(t)
 \geq 0  
\end{array}
\\[6pt]
\displaystyle \int_0^\infty \ga_s(\tau) d \tau = \int_0^{t_0} \ga_s(\tau) d \tau = s \in
[0,2\pi]\, .
\end{array}
$$
The vector field $\bu$ from (\ref{3.1}) is divergence free and
smooth (in Cartesian coordinates $(x,y,z)$). Evidently, there exist $\bff^E$ 
and $\bff^N$, smooth axially symmetric vector fields such that $\bu$ satisfies
(with constant pressure) the Euler equations and the Navier-Stokes equations,
respectively. 

In the cylinder $|z|\leq Z_i$, $r\leq R_i$
it corresponds to the rotation by the angle $s$ during the time interval $[0,t_0]$
and outside of the cylinder $|z|\leq Z^o$, $r\leq R^o$ the fluid does not
move at all.

Let us start by analyzing $\bA^E$.  This is actually quite easy to
compute explicitly.
Writing the input vector in cylindrical coordinates, and the output
in Cartesian coordinates, we have
$$
\bA^E(r,z,\vt, t_0) = \Big(r \cos [\vt - s \al (r)\be (|z|)], 
r \sin[\vt - s \al (r) \be (|z|) ], z\Big)
$$
that is,
%\begin{equation} \label{3.2}
$$
\begin{array}{lcr}
\bA^E(x,y,z,t_0)& = &\Big( x \cos[s \al (\sqrt{x^2+ y^2})\be (|z|)]
 +  y \sin[s \al (\sqrt{x^2+ y^2})\be (|z|) ], \\[4pt]
& & -x \sin[s \al (\sqrt{x^2+ y^2})\be (|z|)] + y \cos[s \al (\sqrt{x^2+ y^2})
\be (|z|) ] , z \Big)\, .
\end{array}
%\end{equation}
$$
Inside the inner cylinder we have
$$
\nabla \bA^E =
\left( \begin{array}{ccc}
\cos s, & \sin s, & 0 \\
-\sin s, & \cos s, & 0 \\
0, & 0, & 1 \end{array} \right)\, ;
$$
outside the outer cylinder
$$
\nabla \bA^E =
\left( \begin{array}{ccc}
1, & 0, & 0 \\
0, & 1, & 0 \\
0, & 0, & 1 \end{array} \right)\, ;
$$
for $|z|\leq Z_i$, $R_i \leq r\leq R^o$
$$
\nabla \bA^E =
\left( \begin{array}{ccc}
\cos [s \al (r) ], & \sin [s \al (r) ],  & 0 \\
-\sin [s \al (r)], & \cos [s \al (r) ] , & 0 \\
0, & 0,& 1 \end{array} \right) + \bM_1 +\bM_2
$$
with 
$$
\begin{array}{l}
\bM_1 =
\left( \begin{array}{ccc}
- s\frac{x^2}{r} \sin [s\al (r)]\al'(r), &
s\frac{y^2}{r} \cos [s\al (r)]\al'(r), & 0 \\
-s\frac{x^2}{r} \cos [s\al (r)]\al'(r), & - s\frac{y^2}{r} 
\sin [s\al (r)]\al'(r), & 0 \\
0, & 0, & 0 \end{array} \right)\, ,
\\[6pt]
\bM_2 =
\left( \begin{array}{ccc}
s\frac{xy}{r} \cos [s\al (r)]\al'(r), & - s\frac{xy}{r} 
\sin [s\al (r)]\al'(r), & 0 \\
- s\frac{xy}{r} \sin [s\al (r)]\al'(r), & 
-s\frac{xy}{r} \cos [s\al (r)]\al'(r), & 0 \\
0, & 0, & 0 \end{array} \right)\, ;
\end{array}
$$
for $Z_i\leq |z|\leq Z^o$, $r\leq R_i$
$$
\nabla \bA^E =
\left( \begin{array}{ccc}
1, & 0, & -s x  \, \sign(z) \sin [s \be (|z|)] \be'(|z|)  +s y \, \sign(z) \cos 
[s \be (|z|)] \be'(|z|)\\
0, & 1, &  -s y \, \sign(z) \sin [s \be (|z|)] \be'(|z|)  -sx \, \sign(z) \cos 
[s \be (|z|)] \be'(|z|)\\
0, & 0, & 1 \end{array} \right)\, ;
$$
and finally for $Z_i\leq |z|\leq Z^o$, $R_i \leq r\leq R^o$ we get a combination
of the last two cases. We will use the structure of $\nabla \bA^E$ later.

\bigskip

Let us now look at the difference between $\bA^N$ and $\bA^E$, our goal
being inequality~(\ref{3.4}) below.  We have
$$
\begin{array}{c}
\displaystyle \pder{}{t} (\bA^N - \bA^E) + \bu \cdot \nabla (\bA^N - \bA^E) = \nu \Delta \bA^N
\\ [6pt]
(\bA^N - \bA^{E})(\bx, 0) = \bnul\, .
\end{array}
$$
Taking the spatial gradient we get
\begin{equation} \label{3.3}
\pder{}{t} [\nabla (\bA^N - \bA^E)] + \bu \cdot \nabla [\nabla (\bA^N - \bA^E)]
 = \nu \Delta \nabla \bA^N - (\nabla (\bA^N-\bA^E))\nabla \bu\, .
\end{equation}
Now, since 
$$
\sup _{t \in [0,t_0]}\|\nabla ^3 \bA^N\|_p \leq C (\|\nabla \bu\|_{k,p})
$$
for some $k$ sufficiently large, we have, after testing equation (\ref{3.3})
by $|\nabla (\bA^N-\bA^E)|^{p-2}\nabla (\bA^N - \bA^E)$
$$
\frac{d}{dt} \|\nabla (\bA^N-\bA^E)\|_p \leq \nu \|\nabla^3 \bA^N\|_p +
\|\nabla \bu\|_\infty  \|\nabla (\bA^N-\bA^E)\|_p\, .
$$
Thus, as $\nabla (\bA^N-\bA^E)(\bx,0) = \bnul$, we get
\begin{equation} \label {3.4}
\sup _{t \in [0,t_0]} \|\nabla (\bA^N-\bA^E)\|_p \leq \nu 
C(\|\nabla \bu\|_{k,2}, t_0)
\end{equation}
for all $p\in (1,\infty]$.

\section{The decay of $\nabla \bA^N-\bI$}
\setcounter{equation}{0}

Let us now look in particular at $\nabla \bA^N$ for $t > t_0$. We have that 
$\nabla \bA^N$ satisfies the heat equation
$$
\begin{array}{c}
\displaystyle \pder{}{t} (\nabla \bA^N) - \nu \Delta (\nabla \bA^N) = \bnul
 \qquad {\rm in \ }
\rdd^3 \times (t_0,\infty) \\
\nabla \bA^N(t_0) \quad {\rm given.}
\end{array}
$$

Therefore also
$$
%\begin{array}{c}
\displaystyle \pder{}{t} (\nabla \bA^N - \bI) - \nu \Delta (\nabla \bA^N-\bI) 
= \bnul \qquad {\rm in \ }
\rdd^3 \times (t_0,\infty) 
%\\
%\nabla \bA^N(t_0) \quad {\rm given.}
%\end{array}
$$
and, especially, at $\bx=\bnul$,
\begin{equation} \label{4.1}
\nabla \bA^N(\bnul, t) - \bI = \frac {C}{(t-t_0)^{\frac 32} \nu^{\frac 32}} 
\int_{\rtt^3} \ee^{-\frac{|\bp|^2}{4\nu (t-t_0)}} (\nabla \bA^N-\bI)(\bp,t_0) 
d\bp\, .
\end{equation}

Our first goal will be to show that the function $F$ in Lemma~\ref{homotopy}
satisfies $F(t,s) \to (1,0)$ as $t \to \infty$ uniformly in $s \in [0,2\pi]$.
This will complete the proof that $F$ is continuous on 
$[0,\infty]\times[0,2\pi]$, and that $F(\infty,s) = (1,0)$, so that
$F$ is indeed a homotopy.

We have that $\nabla\bA^N - \bI = (\nabla\bA^N- \nabla\bA^E) +  (\nabla\bA^E
-\bI)$.
Using the fact that $\nabla \bA^E(\bp, t_0) -\bI$ has bounded support and 
$\|\nabla (\bA^N -\bA^E)\|_p(t_0) \leq C$,
we use the H\"older inequality and end up with
$$
|\nabla \bA^N(0, t) - \bI| \leq \frac{C}{(t-t_0)^a}
$$
with some positive power $a$ ($C$ may depend on any constants which
appeared above, but is independent of the time).

Since it is clear that $F(\cdot,0)$ is the constant function,
the proof of Lemma~\ref{homotopy} will be complete when we have
shown that 
$F(\cdot,2\pi)$ is homotopic to the function whose winding number is $1$,
at least if $\nu$, $R^o-R_i$, and $Z^o-Z_i$ are small enough.

It is clear that the representation of 
$\nabla\bA^E_{2\pi}(\cdot,0)$ has this property.  
So let us put $s=2\pi$.  
We need to
show that $\nabla\bA^E_{2\pi}(\cdot,0) - \nabla\bA^N_{2\pi}(\cdot,0)$ is 
small enough
to construct a linear homotopy between the representation of
$\bA^E_{2\pi}(\cdot,0)$ and 
$F(\cdot,2\pi)$ that does not pass through $(0,0)$.  We have already shown
this property for $t \le t_0$ in equation~(\ref{3.4}), at least when $\nu$
is sufficiently small.  So all that remains is to show the following result.

\begin{lema} \label{l 4.1}
There exist $\eps_1$ and $\eps_2 >0$ such that if $\max \{R^o-R_i, 
Z^o-Z_i\} \leq \eps_1$ and $\nu \leq \eps_2(\eps_1)$ then for $s=2\pi$,
$ \nabla \bA^N - \bI)(0,t) \leq \frac{1}{10}$ for any 
$t\geq t_0$.
\end{lema}

\demo
We denote by $I_1$ the part of integral (\ref{4.1})
with  $(\nabla \bA^N-\bI)(\bp,t_0)$ replaced by $(\nabla \bA^N-\nabla
\bA^E)(\bp,t_0)$, and by $I_2$--$I_5$ the parts of integral (\ref{4.1})
with $(\nabla \bA^N-\bI)(\bp,t_0)$ replaced by $(\nabla \bA^E-\bI)(\bp,t_0)$;
namely by $I_2$ the integral over the inner cylinder, by $I_3$ over the cylinder 
$C(R^o,Z_i)$ without the inner cylinder, by $I_4$ the integral over the 
outer cylinder $C(R^o,Z^o)$ minus the cylinder $C(R^o,Z_i)$ and finally by
$I_5$ over the complement of the outer cylinder.

Evidently, $I_2 = I_5 = 0$ since $s=2\pi$. Let us now consider $I_3$.
If we rewrite $\nabla A^E(\bnul, t_0)-\bI$ (in Cartesian components) 
into the cylindrical coordinates, we get that it is equal to $\bM_0 +
\bM_1 + \bM_2$ with
$$
\begin{array}{l}
\bM_0 =
\left( \begin{array}{ccc}
\cos [2\pi\al (r)]-1, &
\sin[2\pi \al (r)], & 0 \\
-\sin [2\pi\al (r)], &  
\cos [2\pi\al (r)]-1, & 0 \\
0, & 0, & 0 \end{array} \right)\, ,
\\[6pt]
\bM_1 =
\left( \begin{array}{ccc}
- 2\pi r \sin [2\pi\al (r)]\al'(r)\cos^2\vt, &
2\pi r  \cos [2\pi \al (r)]\al'(r) \sin^2\vt, & 0 \\
- 2\pi r \cos [2\pi\al (r)]\al'(r) \cos^2\vt, & - 2\pi r 
\sin [2\pi\al (r)]\al'(r) \sin^2\vt, & 0 \\
0, & 0, & 0 \end{array} \right)\, ,
\\[6pt]
\bM_2 =
\left( \begin{array}{ccc}
2\pi r  \cos [2\pi\al (r)]\al'(r) \sin \vt\cos\vt, & - 2\pi r 
\sin [2\pi\al (r)]\al'(r)\sin\vt\cos\vt, & 0 \\
- 2\pi r \sin [2\pi\al (r)]\al'(r)\sin \vt \cos\vt, &  
-2\pi r \cos [2\pi\al (r)]\al'(r) \sin\vt\cos\vt, & 0 \\
0, & 0, & 0 \end{array} \right)\, .
\end{array}
$$

The heat kernel is independent of the angle $\vt$; after integration over 
it the matrix $\bM_2$ disappears and from $\bM_1$ 
we are left with integrals of the type
$$
\frac{C}{(t-t_0)^{\frac 32} \nu^{\frac 32}} \int_{\index {\begin{array}{c}{
 R_i\leq r
\leq R^o }\\ { |z|\leq Z_i}\end{array}}} \ee^{-\frac{r^2+z^2}{4\nu (t-t_0)}} r^2
\sin [2\pi \al(r)] \al'(r) dr dz
$$
(in some terms, $\sin$ is replaced by $\cos$). Using the standard change of 
variables and integrating over the $z$ variable we end up with 
$$
C \int_{\frac{R_i}{\sqrt{\nu(t-t_0)}} \leq u \leq 
\frac{R^o}{\sqrt{\nu(t-t_0)}} } \ee^{-\frac{u^2}{4}} u^2 \sin \Big(2\pi 
\al(u \sqrt{\nu(t-t_0)})\Big) \al'(u\sqrt {\nu(t-t_0)}) d u\, .
$$

Now the application of the Taylor theorem on the function $\ee^{-\frac{u^2}{4}} u^2$
yields 
$$
\ee^{-\frac{u^2}{4}} u^2 = \ee^{- \frac{R_i^2}{4 \nu (t-t_0)}} \frac{R_i^2}
{\nu (t-t_0)} + \ee^{-\frac {\xi^2}{4}} \Big(2\xi - \frac{\xi^3}{2}\Big) 
\Big(u-\frac
{R_i}
{\sqrt{\nu (t-t_0)}}\Big)\, ,
$$
where $\xi \in (\frac {R_i} {\sqrt{\nu (t-t_0)}}, \frac {R^o} 
{\sqrt{\nu (t-t_0)}})$. Moreover
$$
\int_{\frac{R_i}{\sqrt{\nu(t-t_0)}} \leq u \leq 
\frac{R^o}{\sqrt{\nu(t-t_0)}} } \sin \Big(2\pi \al(u \sqrt
{\nu(t-t_0)})\Big) \al'(u\sqrt {\nu(t-t_0)}) d u = 0\, .
$$
A similar argument can be applied also on terms coming from $\bM_0$. Thus
we have
$$
\begin{array}{c}
\displaystyle 
|I_3| \leq C\ee^{-\frac{R_i^2}{4\nu (t-t_0)}} \Big(\frac {2 R^o}{\sqrt{\nu 
(t-t_0)}} + \frac {(R^o)^3}{\Big(\sqrt{\nu 
(t-t_0)}\Big)^3}\Big) \\[6pt]
\displaystyle \max_{u \in (\frac {R_i}{\sqrt{\nu (t-t_0)}},\frac 
{R^o}{\sqrt{\nu (t-t_0)}})}|\al'(u \sqrt{\nu (t-t_0)})| \frac{(R^o-R_i)^2}
{\nu(t-t_0)} +\ee^{-\frac{R_i^2}{4\nu (t-t_0)}} (R^o-R_i)\frac{R^o}
{\sqrt{\nu(t-t_0)}}\, .
\end{array}
$$
We can choose $\al(r)$ in such a way that $\al'(r) \leq \frac{C}{R^o-R_i}$ and
as 
$$
\ee^{-\frac{R_i^2}{4\nu (t-t_0)}} \Big(\sqrt{\nu (t-t_0)}\Big)^a \leq C(a, R_i)
$$
for any $a \in \rdd$, we finally get
$$
|I_3| \leq C (R^o-R_i)
$$
with the constant in particular independent of $\nu$ and $t$. Therefore for 
the "boundary layer" sufficiently thin, this term can be done arbitrarily 
small, independently of the viscosity and the time.

Similarly we can estimate $I_4$; here $(\nabla \bA^E)_{i,3}^{i=1,2}(\bp,t_0)$
are odd functions in $z$ and thus we get zero after the integration of 
the $z$ variable. 
For the components $i,j$; $i,j=1,2$ proceed similarly as above and end up with 
the following integral
$$
\int_\Omega \ee^{-\frac{u^2+v^2}{4}} 
u^2 \sin \Bigg(2\pi \al\Big[u \sqrt
{\nu(t-t_0)}\Big] \be\Big[|v|\sqrt{\nu (t-t_0)\Big]}\Bigg) \al'\Big(u\sqrt {\nu(t-t_0)}
\Big) d u d v
$$
with $\Omega =\{(u,v); \frac{R_i}{\sqrt{\nu(t-t_0)}} \leq u \leq 
\frac{R^o}{\sqrt{\nu(t-t_0)}},   \frac{Z_i}{\sqrt{\nu (t-t_0)}} \leq |v|
\leq \frac{Z^o}{\sqrt{\nu(t-t_0)}}\}$. 
But now $\be \neq 1$ and we cannot proceed as above. Nevertheless, we get
that the integral above is bounded by
$$
C (Z^o-Z_i) (\nu (t-t_0))^a \ee^{-\frac{R_i^2 + Z_i^2}{4\nu (t-t_0)}}
$$
with the constant independent of $\nu$ and $t$.
Thus, if $Z^o-Z_i$ is small, we get that also $I_4$ is small.

Finally,
$$
\begin{array}{c}
\displaystyle |I_1|  \leq \frac{C}{\Big(\nu(t-t_0)\Big)^{\frac 32}} \int_{\rtt^3} \ee^{-\frac
{|\bp|^2}{4\nu (t-t_0)}} |\nabla \bA^N - \nabla \bA^E|(\bp,t_0) d\bp \\[6pt]
\leq
C \|\nabla \bA^N - \nabla \bA^E\|_\infty(t_0) \leq \nu C(R^o-R_i, Z^o-Z_i, t_0)\,.
\end{array}
$$
\kondemo

\begin{rema} \label {r 5.2}
{\rm We have shown that $\bA^N$ may have no smooth inverse.  However it would be
more interesting to provide an example in which it can be shown that
$\bA^N$ has no inverse at all.
Looking at the representation of $\nabla \bA^E_s(\bnul, t_0)$ it is not 
difficult to see 
that $\nabla^2 \bA_s^E(\bnul, t_0)$ is odd in $x$ and $y$ and therefore, 
since
the same holds also for $\nabla^2 \bA_s^N(\bnul, t_0)$, we get
that  $\nabla^2 \bA_s^N(\bnul, t) = \bnul$ for any $t>t_0$ and any $s \in 
[0,2\pi]$ and thus
$\bA^N$ is in fact invertible with a non-smooth inverse.}
\end{rema}

\setcounter{equation}{0}
\section{The Probabilistic Approach}

Here we will describe a probabilistic approach to solving 
equation~(\ref{navier-stokes-m}).  For simplicity let us
consider the case when the forcing term $\bff = 0$.  
We will not be rigorous.

We will
suppose that we have found $\bu$ using equation~(\ref{navier-stokes-u}).
Now let $\bb_t$ be a Brownian motion in 3 dimensions, starting at
the origin.  Define $\tilde\bu(\bx,t) = \bu(\bx + 2\nu \bb_t,t)$.
Let $\widetilde \bm$ be a random vector field that satisfies
the equations
\begin{equation} \label{euler-tilde-m} 
\begin{array}{c}
\displaystyle \pder{\widetilde\bm}{t} + \tilde\bu\cdot \nabla \widetilde\bm 
  + \widetilde\bm \cdot (\nabla \tilde\bu)^T = 
  0
\quad {\rm in \ } \rdd^3 \times (0,T) \\[10pt]
\widetilde \bm(\bx,0) = \bu_0(\bx) \quad {\rm in \ } \rdd^3\, .
\end{array}
\end{equation}
Now let $\overline\bm(\bx,t) = \bm(x-2\nu\bb_t,t)$.  Then
$\bm(\bx,t) = E(\overline\bm(\bx,t))$ satisfies 
equation~(\ref{navier-stokes-m}).
(Here $E(\cdot)$ represents the expected value.)

The reason why this works is because of the It\^o formula.  We have
that
$$ \pder{\overline\bm}{t} + \bu\cdot \nabla \overline\bm 
   + \overline\bm \cdot (\nabla \bu)^T = 
   \nu \Delta \overline \bm + 2 \nu \pder{\bb}t \cdot \nabla \overline\bm ,$$
and taking expectations the result follows.

The solution to equation~(\ref{euler-tilde-m}) can be computed as 
follows.  Suppose that the initial value of $\bm$ satisfies
equation~(\ref{init-m-alpha-beta}).  Then if 
$$
\begin{array}{c}
\displaystyle \pder{\alpha_i}{t} + \tilde\bu \cdot \nabla \alpha_i = 0
\\[6pt]\\
\displaystyle \pder{\beta_i}{t} + \tilde\bu \cdot \nabla \beta_i =0
\end{array}
$$
then 
$$
\widetilde \bm(\bx,t) = \sum_{i=1}^R \beta_i(\bx,t) \nabla \alpha_i(\bx,t) 
$$
is the solution to system~(\ref{euler-tilde-m}).  But the transport
equations are easily solved by
$\alpha_i(\bx,t) = \alpha_i(\tilde\bA(\bx,t),0)$ and
$\beta_i(\bx,t) = \beta_i(\tilde\bA(\bx,t),0)$, where $\tilde\bA$ is the
back to coordinates map induced by the flow $\tilde \bu$.

This can be used to obtain the following plausibility argument for the
regularity of the Navier-Stokes equations.  Let $W^{-1,BMO}$ denote the
space of functions from $\rdd^3$ for which minus one derivative is in
the space of functions of bounded mean oscillation.  It is known that the
space $L^{\infty}(I; W^{-1,BMO})$ is a critical space for proving
regularity for the Navier-Stokes
equations (see below).  That is, if one can show that the solution to the Navier-Stokes
equations is uniformly in time in any space better than $W^{-1,BMO}$ (such as $W^{-1+\epsilon,BMO}$
for any $\epsilon>0$), then the solution is regular.

Now if the initial data are very nice, then by using some partition of unity
argument, we may suppose that indeed the initial value of $\bm$ does satisfy
equation~(\ref{init-m-alpha-beta}) for some finite value of $R$, 
where the initial values of  $\alpha_i$ and $\beta_i$ are compactly
supported smooth functions.  Then it is easy to see that the solutions
for $\alpha_i$ and $\beta_i$ provided by the transport equations stay
uniformly in $L^\infty$.  Thus it follows that $\nabla \alpha_i$ is uniformly
in the space $W^{-1,BMO}$.

Thus $\widetilde \bm$ is a finite sum of a product of functions uniformly in $L^\infty$
and functions uniformly in $W^{-1,BMO}$.  Thus it might seem that we
are close to showing that $\bu$ (which is the Leray projection of an average
of translations of $\tilde \bm$) is in a space that is critical for proving
regularity.

There are some large, probably insurmountable problems with this approach.
The lesser problem is that we need a space that is better than critical.
The bigger problem is that the space created by taking the convex closure
of products of bounded
functions and functions in $W^{-1,BMO}$ is not really a well defined
space, in that it encompasses every function.

{\bf Criticality of $L^{\infty}(I; W^{-1,BMO})$:} Let us present a
formal proof of this fact, in the case of the Cauchy problem with zero
right-hand side.
Let $\bu$ be the solution to the Navier-Stokes
equations which belongs to the space $L^{\infty}(I;
W^{-1,BMO})$. Multiply equation (\ref{navier-stokes-u})$_1$ by $\Delta
\bu$ and integrate over $\rdd^3$. Notice also that
$$
\left| \int \Delta \bu \cdot (\bu \cdot \nabla \bu) \right| =
\left|\int \frac\partial{\partial x_i} \bu \cdot 
 \left(\frac\partial{\partial x_i} \bu \cdot \nabla \right)\bu\right|
\le \| \nabla \bu \|_3^3 . $$
Then
$$
\frac 12 \frac{d}{dt} \|\nabla \bu\|_2^2 + \nu \|\nabla ^2 \bu\|_2^2
\leq \| \nabla \bu \|_3^3 .
$$
Using the inequality
$$
\|\nabla \bu\|_3 \leq C \|\bu\|_{-1,BMO}^{\frac 13} \|\nabla^2 \bu\|_2^{\frac 23}
$$
(see \cite{Oru}) we get
\begin{equation}
\label{interp}
 \frac{d}{dt} \|\nabla \bu\|_2^2 + \nu \|\nabla ^2 \bu\|_2^2
 \leq
 C \|\bu\|_{-1,BMO} \|\nabla^2 \bu\|_2^2     
\end{equation}
and if $\|\bu\|_{-1,BMO}$ is sufficiently small, the solution is
smooth.
The proof can be done rigorously using the fact that for smooth
initial condition there exists a local smooth solution; on this
interval we obtain estimate (\ref{interp}) and therefore the solution
can not blow up.

\noindent
{\it Stephen Montgomery-Smith, University of Missouri,
Department of Mathematics,  Math. Science Building, 65203 Columbia, MO,
USA} \newline (e-mail {\tt stephen@math.missouri.edu}
\newline
URL {\tt http://www.math.missouri.edu/\~{}stephen}).

\medskip
\noindent
{\it Milan Pokorn\'{y}, Math. Institute of Charles University, 
Sokolovsk\'a 83, 186 75 Praha 8, Czech Republic} \newline (e-mail {\tt
pokorny@karlin.mff.cuni.cz}).

\end{document}